\documentclass[a4paper,10pt]{amsart}
\usepackage[utf8]{inputenc}
\usepackage{amsmath,amssymb,amsthm}
\usepackage{xcolor}
\usepackage{enumerate}
\usepackage{tikz}
\usepackage{tikz-cd}
\usepackage[all]{xy}
\usepackage{float}
\usepackage{soul}
\usetikzlibrary{decorations.pathreplacing}

\usepackage{listings}

\usepackage{hyperref}
\hypersetup{hidelinks}

\usepackage[noabbrev,capitalise,nameinlink]{cleveref}
\crefname{subsection}{Subsection}{Subsections}

\usepackage{thmtools}

\usepackage{stmaryrd}

\theoremstyle{plain}

\newtheorem{thm}{Theorem}[section]

\newtheorem*{conjecture*}{Conjecture}
\DeclareMathOperator{\Aut}{Aut}
\DeclareMathOperator{\supp}{supp}
\DeclareMathOperator{\pr}{pr}
\DeclareMathOperator{\stab}{Stab}
\newcommand{\Pins}{P_{\mathrm{ins}}}
\newcommand{\Prob}{\mathbb{P}}
\newcommand{\Expect}{\mathbb{E}}

\theoremstyle{definition}

\title{Asymptotic insolubility of alternating groups}

\author{Andrea Lucchini}
\address{Andrea Lucchini. University of Padova (Italy), Dipartimento di Matematica ``Tullio Levi Civita''. ORCID: https://orcid.org/0000-0002-2134-4991}
\email{lucchini@math.unipd.it}

\date{}

\begin{document}
	\begin{abstract}
	Let $a,b\in\mathrm{Aut}(A_n)$ with $a\neq 1$, and let $s$ be chosen
	uniformly in $A_n$. We show that the probability that $\langle a, bs\rangle$
	is insoluble tends to $1$ as $n\to\infty$, uniformly in $a$ and $b$. More
	precisely, this probability is at least $1-O(n^{-1/24})$. The proof combines
	a bound in terms of the support of $a$ with estimates on permutations lying
	in transitive subgroups of the symmetric group.	
	\end{abstract}
\maketitle

\section{Introduction}
The following result is stated in \cite{fgg}: {\sl{
There exists an absolute constant $\epsilon > 0$ such that the following holds. 
Let $S$ be a nonabelian finite simple group and let $x, y \in \Aut(S)$ with $x\neq 1$. Then the 
probability that $x$ and a random element of the coset $Sy$ generate a nonsoluble group is at least $\epsilon.$}} The statement was proved in \cite{fgg} for simple groups of Lie type, while the alternating case was considered in \cite{Lucchini}. More precisely, 
for \(n\ge 5\) and \(a,b\in\Aut(A_n)\), let \(\Pins(A_n,a,b)\) be the probability
 that \(\langle a,bs\rangle\) is insoluble when \(s\) is chosen uniformly
 from the inner automorphism group \(A_n\). 
 \cite[Theorem~6]{Lucchini} gives a positive lower bound on this
 probability, uniform in \(n\), \(a\ne 1\) and \(b\). In this note we revise \cite[Theorem~6]{Lucchini} giving a more effective proof and a stronger statement.
 
 \begin{thm}\label{main}
 	Set
 	\[
 	I_n=\min_{\substack{a,b\in\Aut(A_n)\\a\ne 1}}\Pins(A_n,a,b).
 	\]
 	Then \(I_n\to 1\) as \(n\to\infty\). More precisely,
 	\[
 	1-I_n=O(n^{-1/24}).
 	\]
 \end{thm}
 
\section{Proof of \cref{main}}
 	For \(n\ge 7\), identify \(\Aut(A_n)\) with \(S_n\) via conjugation.
 	Fix \(1\ne\sigma\in S_n\), let \(\Sigma=\supp(\sigma)\) and \(m=|\Sigma|\),
 	and choose \(\tau\) uniformly from \(S_n\). Write
 	\(H_\tau=\langle\sigma,\tau\rangle\).
 	
 	\medskip
 	\noindent\emph{A bound in terms of the support.}
 	Fix $x\in \Sigma$. Let $K$ be a subgroup of $H_\tau$ minimal with respect to the property that $\langle \stab_{H_\tau}(x), \sigma\rangle \leq K$. We consider the action of $K$ on the $K$-orbit $\Delta=x^K$ of $x$. Let $J$ be a maximal subgroup of $K$ containing $H_x$. If $T={t_1,\dots,t_r}$ is a right transversal of $J$ in $K$, then $\mathcal B={x^{Jt_1},\dots,x^{Jt_r}}$ is a system of blocks for the action of $K$ on $\Delta$, which is permuted primitively by $K$. Every nonidentity element of a soluble primitive permutation group moves at least half of the points. Since $\sigma\notin J$, it follows that $\sigma$ moves at least half of the blocks in $\mathcal B$, and consequently it also moves at least half of the elements of $\Delta$ (Indeed, a minimal normal subgroup of
 	such a group is elementary abelian and regular, so the group acts by
 	affine transformations. The fixed points of a nonidentity affine
 	transformation, if any, lie in a proper affine subspace). This implies $|\Delta|\leq 2m.$

 	Let \(t\) be the least positive integer with \(\tau^t\in K\).
 	If \(x\) and \(z=x^\sigma\) lie in the same cycle of \(\tau\), write \(u\)
 	and \(v\) for the forward cyclic distances from \(x\) to \(z\) and
 	from \(z\) to \(x\), respectively. Then
 	\(\tau^u\sigma^{-1}\in \stab_{H_\tau}(x)\) and \(\tau^{u+v}\in \stab_{H_\tau}(x)\), so \(t\)
 	divides \(u\) and \(v\). The \(\tau^t\)-orbit of \(x\) contributes
 	\((u+v)/t\) points to \(\Delta\).
 	If \(x\) and \(z\) lie in different cycles, let \(u\) and \(v\)
 	be their respective cycle lengths. Now \(\tau^u\in \stab_{H_\tau}(x)\) and
 	\(\tau^v\in \stab_{H_\tau}(z)=\sigma^{-1}\stab_{H_\tau}(x)\sigma\), and both stabilizers are
 	contained in \(K\). Hence \(t\) again divides \(u\) and \(v\), and
 	the two disjoint \(\tau^t\)-orbits contribute \((u+v)/t\) points
 	to \(\Delta\). In either case, \(s=(u+v)/t\le 2m\).
 	
 	Fix positive integers \(u,v\) with \(l=u+v\le n\).
 	If \(x\) and \(z\) lie in the same cycle of \(\tau\), write the \(l\)-cycle starting at \(x\) and
 	place \(z\) exactly \(u\) steps later. In the other case, write the
 	\(u\)-cycle starting at \(x\) and the \(v\)-cycle starting at \(z\).
 	Both constructions leave \(l-2\) ordered positions to fill.
 	There are \((n-2)!/(n-l)!\) ways to fill them with distinct points,
 	and \((n-l)!\) ways to permute the unused points. Thus each case
 	gives exactly \((n-2)!\) permutations for the prescribed pair,
 	with probability \(1/[n(n-1)]\).
 	
 	For a fixed \(s=(u+v)/t\), there are \(s-1\) choices for \(u/t\),
 	which determines \(v/t\), and at most \(\lfloor n/s\rfloor\) choices
 	for \(t\), since \(ts\le n\). This gives at most \(n\) pairs for each
 	\(s\) in each case. Summing over \(2\le s\le 2m\) gives at most
 	\(2m/(n-1)\) for each case. Adding the two contributions gives
 	\begin{equation}
 		\label{eq:support-bound}
 		\Prob(H_\tau\text{ is soluble})\le \frac{4m}{n-1}.
 	\end{equation}

 	\medskip
 	\noindent\emph{A bound from the orbit meeting the support.}
 	Let \(\Omega_\tau\) be the union of the cycles of \(\tau\) that meet \(\Sigma\).
 	The projection  $\pr_\Sigma: S_n\to \operatorname{Sym}(\Sigma)\) which maps $\tau$ to
 	the permutation \(\tau_\Sigma\), obtained by deleting the points outside
 	\(\Sigma\) from the cycles of \(\tau\), is uniform on \(\operatorname{Sym}(\Sigma)\), i.e. the size of $\pr_\Sigma^{-1}(\rho)$ is the same for all $\rho \in \operatorname{Sym}(\Sigma)$ \cite[Observations~11]{BabaiHayes}.
 	The orbits of \(\langle\sigma|_\Sigma,\tau_\Sigma\rangle\) are the nonempty
 	intersections of \(\Sigma\) with the \(H_\tau\)-orbits
 	\cite[Observations~12]{BabaiHayes}. Since \(\sigma|_\Sigma\) has no
 	fixed points, the transitivity bound of Babai and Hayes
 	\cite[Corollary~9]{BabaiHayes} shows that \(\Omega_\tau\) is a single \(H_\tau\)-orbit
 	except with probability \(O(1/m)\).
 	
 	For \(w\notin \Sigma\), the cycle of \(w\) under \(\tau\) misses \(\Sigma\)
 	exactly when \(w\) is fixed by the projection of \(\tau\) onto
 	\(\Sigma\cup\{w\}\). This projection is uniform, so
 	\[
 	\Expect(n-|\Omega_\tau|)=\frac{n-m}{m+1}.
 	\]
 	Applying Markov's inequality to \(n-|\Omega_\tau|\) gives
 	\begin{equation}
 	\Prob(|\Omega_\tau|<n/2)
 	=\Prob(n-|\Omega_\tau|>n/2)
 	\le \frac{2}{n}\Expect(n-|\Omega_\tau|)
 	\le \frac{2}{m+1}.
 	\end{equation}
 	
 	For a set \(U\) of size \(k\), let \(\mathcal T_U\) be the set of
 	permutations lying in a transitive subgroup of \(\operatorname{Sym}(U)\)
 	other than \(\operatorname{Alt}(U)\) or \(\operatorname{Sym}(U)\), and set
 	\(T(k)=|\mathcal T_U|/k!\).
 	If \(H_\tau\) is soluble and \(\Omega=U\) is a single \(H_\tau\)-orbit of size \(k\ge5\),
 	then \(\tau|_U\in\mathcal T_U\). Indeed, the group induced by \(H_\tau\) on
 	\(U\) is soluble and transitive. Both \(\operatorname{Alt}(U)\)
 	and \(\operatorname{Sym}(U)\) are insoluble when \(k\ge5\).
 	We claim that
 	for every fixed \(U\supseteq \Sigma\) of size \(k\ge5\), we have
 	\[
 	\Prob\bigl(\tau|_U\in\mathcal T_U\mid \Omega_\tau=U\bigr)
 	\le \frac{k}{m}T(k).
 	\]
 	To prove this bound, let \(\mathcal C_U\) consist of the permutations
 	of \(U\) whose cycles all meet \(\Sigma\). The condition \(\Omega_\tau=U\) holds
 	exactly when \(\tau\) preserves \(U\) and \(\tau|_U\in\mathcal C_U\).
 	Each member of \(\mathcal C_U\) extends to exactly \((n-k)!\)
 	permutations with \(\Omega_\tau=U\), one for each permutation of the
 	complement of \(U\). Counting these extensions gives
 	\[
 	\Prob\bigl(\tau|_U\in\mathcal T_U\mid \Omega_\tau=U\bigr)
 	=\frac{|\mathcal T_U\cap\mathcal C_U|}{|\mathcal C_U|}
 	\le\frac{|\mathcal T_U|/k!}{|\mathcal C_U|/k!}.
 	\]
 	Here \(|\mathcal C_U|/k!=m/k\) by \cite[Lemma~10]{Lucchini}.

 	For \(n\ge10\), average over the possible values \(U\) of \(\Omega_\tau\)
 	with \(|U|\ge n/2\). Since \(k/m\le n/m\), the probability that \(H_\tau\)
 	is soluble and \(\Omega_\tau\) is a single \(H_\tau\)-orbit of size at least \(n/2\)
 	is at most \((n/m)\sup_{k\ge n/2}T(k)\). The probability that \(\Omega_\tau\)
 	is not a single orbit or that \(|\Omega_\tau|<n/2\) is \(O(1/m)\). Therefore
 	\begin{equation}
 		\label{eq:orbit-bound}
 		\Prob(H_\tau\text{ is soluble})
 		\le O(1/m)+\frac{n}{m}\sup_{k\ge n/2}T(k).
 	\end{equation}
 	
 	Since $\mathbb{P}(H_\tau \text{ is soluble})$ is bounded above by both the
 	right-hand side of \eqref{eq:support-bound} and that of \eqref{eq:orbit-bound}, it is at most their
 	geometric mean. The factor $m$ cancels, and we obtain
 	\[
 	\mathbb{P}(H_\tau \text{ is soluble})
 	\le \sqrt{\frac{4m}{n-1}\cdot\frac{C+n\sup_{k\ge n/2}T(k)}{m}}
 	= O\!\left(\left(\frac{1}{n}+\sup_{k\ge n/2}T(k)\right)^{1/2}\right),
 	\]
 	uniformly in $\sigma \neq 1$.
 	The quantitative estimate of Eberhard, Ford and Koukoulopoulos
 	\cite[Theorem~1.1]{EFK} gives, in particular, \(T(k)=O(k^{-1/12})\).
 	Finally, conditioning on the parity of \(\tau\) to make it uniform
 	on either coset \(bA_n\) multiplies the probability of solubility by
 	at most \(2\). This gives the asserted bound on \(1-I_n\).

\subsection*{Declaration of generative AI and AI-assisted technologies in the manuscript preparation process}

During the preparation of this work the authors used ChatGPT-6.2 Sol in order to collect and refine ideas and review the manuscript. After using this tools, the authors reviewed and edited the content as needed and take full responsibility for the content of the published article.

\end{document}